\documentstyle[11pt]{article}
\title{Commensurability and locally free Kleinian groups}
\author{James W. Anderson \\ \footnotesize{Faculty of Mathematical
Studies} \\ \footnotesize{University of Southampton} \\
\footnotesize{Southampton SO17 1BJ, England} \\
\footnotesize{jwa@maths.soton.ac.uk}}
\date{1991 Mathematics Subject Classification 57M50, 30F40, 20E25}

\newtheorem{conjecture}{Conjecture}[section]
\newtheorem{corollary}[conjecture]{Corollary}

\newtheorem{proposition}[conjecture]{Proposition}

\newtheorem{remark}[conjecture]{Remark}
\newtheorem{theorem}[conjecture]{Theorem}

\newenvironment{proof}{{\bf Proof} \hspace{.1cm}}{\hfill $\Box$ \\
\vspace{.15cm}} 

\begin{document}

\maketitle

\section{Introduction}
\label{introduction}

The purpose of this note is to show that there exist infinitely many
commensurability classes of finite volume hyperbolic $3$-manifolds whose
fundamental group contains a subgroup which is locally free but not
free.  This gives a strong answer to a question of Kropholler given in
the Problem List in Niblo and Roller \cite{niblo-roller}.

\medskip
\noindent
{\bf Theorem \ref{infinitely-many}: } {\em There exist infinitely many
commensurability classes of finite volume hyperbolic $3$-manifolds
whose fundamental groups contain a subgroup which is locally free but
not free.} 
\medskip

There are three main observations which make up the proof.  The first
observation, discussed in Section \ref{persistence}, is that convex
co-compact subgroups of a fundamental group of a hyperbolic
$3$-manifold $N$ persist in the approximates given by hyperbolic Dehn
surgery. This result is stated formally as Proposition
\ref{persistence-prop}. 

The second observation, discussed in Section
\ref{non-commensurability-criterion}, is that a collection of distinct
hyperbolic $2$- or $3$-manifolds of uniformly bounded volume contains 
infinitely many commensurability classes.   This result is stated
formally as Proposition \ref{commensurability-prop}.  Such collections
of hyperbolic $3$-manifolds are generated by hyperbolic Dehn surgery.

The third observation, discussed in the proof of Theorem
\ref{infinitely-many}, is a general construction of finite volume
hyperbolic $3$-manifolds which is a direct application of Thurston's
hyperbolization theorem for Haken atoroidal $3$-manifolds.  The
material in this Section is well-known, though to my knowledge the
construction given has never been put to paper and so is included for
the sake of completeness. 

\medskip

At this point, given the prevalence of hyperbolic $3$-manifolds whose
fundamental groups contain subgroups which are locally free but not
free, I would like to put forward two complementary conjectural
pictures. One is that the result of Theorem \ref{infinitely-many}
holds not just for infinitely many commensurability classes, but in
fact holds for all finite volume hyperbolic $3$-manifolds.

\begin{conjecture} The fundamental group of every finite volume
hyperbolic $3$-manifold contains a subgroup which is locally free but
not free.
\label{conjecture-all}
\end{conjecture}

Though there is no direct evidence for Conjecture
\ref{conjecture-all}, it is not unplausible, given Theorem
\ref{infinitely-many}, together with Proposition
\ref{commensurable-wacky}, which states that the property of
containing a subgroup which is locally free but not free is an
invariant of a commensurability class.

The second conjecture is based on the as yet unrealized hope that the
property of containing a subgroup which is locally free but not free
can be used to distinguish between certain classes of finite volume
hyperbolic $3$-manifolds.  We note here that Reid \cite{reid-trace}
has shown that there are infinitely many commensurability classes of 
hyperbolic $3$-manifolds which fiber over the circle.

\begin{conjecture} The fundamental group of a finite volume hyperbolic
$3$-manifold $N$ contains a subgroup which is locally free but not
free if and only if $N$ does not have a finite cover which fibers over
the circle.  
\label{conjecture-fiber}
\end{conjecture}

Though there is no direct evidence to support Conjecture
\ref{conjecture-fiber}, one piece of heurstic evidence is that there
are other algebraic properties of the fundamental groups of finite
volume hyperbolic $3$-manifolds which conjecturally fail if and only
if the $3$-manifold in question is commensurable to one that fibers over
the circle.  One such is the {\em finitely generated intersection
property}; for a discussion, see Canary \cite{canary-survey} and the
references contained therein.

Note that if Conjecture \ref{conjecture-fiber} were to hold true, then
all the hyperbolic $3$-manifolds constructed in this note would be
counterexamples to the conjecture of Thurston \cite{thurston-bull}
that a finite volume hyperbolic $3$-manifold has a finite  cover
which fibers over the circle.

\medskip

We close this introduction by giving some definitions.  A {\em
Kleinian group} is a discrete subgroup of ${\rm PSL}_2({\bf C})$,
which we view as acting both on the Riemann sphere $\overline{{\bf
C}}$ by M\"obius transformations  and on real hyperbolic $3$-space
${\bf H}^3$ by isometries, where the two actions are linked by the
Poincar\'e extension.    For the duration of this note, we assume that
{\em all Kleinian groups are torsion-free}.  Note though that by
Selberg's lemma \cite{selberg}, which implies that a finitely
generated Kleinian group contains a torsion-free subgroup of finite
index, that the results of this note hold for Kleinian groups with
torsion as well.

The action of an infinite Kleinian group $\Gamma$ partitions
$\overline{{\bf C}}$ into two sets, the {\em domain of discontinuity}
$\Omega(\Gamma)$, which is the largest open subset of $\overline{{\bf
C}}$  on which $\Gamma$ acts discontinuously, and the {\em limit set}
$\Lambda(\Gamma)$.  If $\Lambda(\Gamma)$ contains two or fewer points,
$\Gamma$ is  {\em elementary}, otherwise $\Gamma$ is {\em
non-elementary}.  For a non-elementary Kleinian group $\Gamma$, the
limit set $\Lambda(\Gamma)$ can also be described as the smallest
non-empty closed subset of $\overline{{\bf C}}$ invariant under
$\Gamma$.  We refer the reader to Maskit \cite{maskit-book} or
Matsuzaki and Taniguchi \cite{matsuzaki-taniguchi} for the basic
elements of the theory of Kleinian groups.

Two Kleinian groups $\Gamma$ and $\Gamma'$ are {\em commensurable} if
their intersection $\Theta =\Gamma\cap\Gamma'$ has finite index in
both $\Gamma$ and $\Gamma'$.  Equivalently, two hyperbolic $3$-manifolds
$N$ and $N'$ are {\em commensurable} if there exists a hyperbolic
$3$-manifold $P$ which is a finite cover of both $N$ and $N'$.  In the
latter case, there are then realizations of both $\pi_1(N)$ and
$\pi_1(N')$ as commensurable Kleinian groups.

\medskip

A finitely generated Kleinian group is {\em convex co-compact} if the
convex core of its associated hyperbolic $3$-manifold ${\bf
H}^3/\Gamma$ is compact.  Recall that the {\em convex core} of a
hyperbolic $3$-manifold $N$ is the minimal convex submanifold of $N$
whose inclusion into $N$ is a homotopy equivalence.   Equivalently, a
Kleinian group $\Gamma$ is convex co-compact if and only if the
associated $3$-manifold $({\bf H}^3\cup\Omega(\Gamma))/\Gamma$ is compact.

More generally, a finitely generated Kleinian group is {\em
geometrically finite} if the convex core of its associated hyperbolic
$3$-manifold has finite volume.  This is one of several equivalent
definitions of geometric finiteness; the interested reader is referred
to Bowditch \cite{bowditch} for a complete discussion.

A Kleinian group $\Gamma$ is {\em minimally parabolic} if every
parabolic element of $\Gamma$ is contained in a ${\bf Z}\oplus {\bf
Z}$ subgroup of $\Gamma$. 

\medskip

An {\em incompressible surface} in a $3$-manifold $M$ is an embedded
orientable surface $S$ in $M$ for which $\pi_1(S)$ is infinite and for
which $\pi_1(S)$ injects into $\pi_1(M)$ under the inclusion map.  If
$M$ is compact, we say that $M$ has {\em incompressible boundary} if
each component of $\partial M$ is an incompressible surface.  A
surface $S$ in $M$ is {\em essential} if it is incompressible and is
not homotopic into $\partial M$.
 
Let $M$ be a compact $3$-manifold and let $S$ be a (possibly
disconnected) compact subsurface of $\partial M$.  Say that $S$ is
{\em an-annular} if there does  not exist an essential annulus
$A$ in $M$ both of whose boundary components are contained in $S$.
A compact $3$-manifold $M$ is {\em acylindrical} if $\partial M$ is
an-annular. 

A compact $3$-manifold $M$ is {\em irreducible} if every embedded
2-sphere in $M$ bounds a $3$-ball in $M$, and $M$ is {\em atoroidal} if
$M$ does not contain an essential torus. 

A compact $3$-manifold $M$ is {\em hyperbolizable} if there exists a
hyperbolic $3$-manifold $N ={\bf H}^3/\Gamma$ homeomorphic to the
interior of $M$.  Note that a hyperbolizable $3$-manifold $M$ is 
necessarily orientable, irreducible, and atoroidal.  The converse
holds for compact $3$-manifolds which contain an incompressible
surface, by the hyperbolization theorem of Thurston, see Morgan
\cite{morgan}.  Also, since the universal cover ${\bf H}^3$ of $N$ is
contractible, the fundamental group of $M$ is isomorphic to
$\Gamma$. For a discussion of the basic theory of $3$-manifolds, we
refer the reader to Hempel \cite{hempel}.

Given a compact, hyperbolizable $3$-manifold $M$, a {\em uniformization}
of $M$ is a Kleinian group $\Gamma$ so that $N ={\bf H}^3/\Gamma$ is
homeomorphic to the interior of $M$. For the purposes of this note, we
restrict our attention to {\em minimally parabolic, geometrically
finite uniformations}.  Note that if $\partial M$ contains no tori and
if $\Gamma$ is a minimally parabolic, geometrically finite
uniformization of $M$, then $\Gamma$ contains no ${\bf Z}\oplus {\bf
Z}$ subgroups. 

\medskip

This paper was written while I was visiting Rice University, and I
like to thank the department there for its hospitality.  I would also
like to thank Graham Niblo, Alan Reid, Joe Masters, and Dick Canary
for helpful conversations about the ideas in this note.

\section{Persistence of convex co-compact subgroups}
\label{persistence}

We begin this Section by describing the operation of hyperbolic Dehn
surgery.  Let $M$ be a compact hyperbolizable $3$-manifold whose
boundary $\partial M$ contains a non-empty collection of tori ${\cal
T} =\{ T_1,\ldots, T_n\}$.   

For each $1\le k\le n$, choose a meridian-longitude system $(m_k,
l_k)$ on $T_k$.  Given a pair $(p_k, q_k)$ of relatively prime
integers, form a new manifold by attaching a solid torus $V$ to $M$ by
an orientation-reversing homeomorphism $g: \partial V \rightarrow T_k$
so that, if $c$ is the meridian of $V$, then $g(c)$ is a $(p_k, q_k)$ 
curve on $T_k$ (with respect to the chosen meridian-longitude
system).  If we write the $n$-tuple of pairs as $({\bf p},{\bf q})
=((p_1, q_1),\ldots, (p_n, q_n))$, we say that $M({\bf p}, {\bf q})$
is obtained from $M$ by {\em $({\bf p}, {\bf q})$-Dehn filling along
$\cal T$}. 

Let $\Gamma$ be a minimally parabolic, geometrically finite Kleinian
group uniformizing $M$, so that the hyperbolic $3$-manifold $N ={\bf
H}^3/\Gamma$ is homeomorphic to the interior of $M$.   Let  $\{ ({\bf
p}_j, {\bf q}_j) = ((p^1_j,q^1_j),  \ldots,(p^n_j,q^n_j)) \}$ be a
sequence of $n$-tuples of pairs of relatively prime integers such
that, for each $1\le k\le n$, $\{ (p^k_j,q^k_j)\}$ converges to
$\infty$ as $j\rightarrow\infty$.   

The generalized hyperbolic Dehn surgery theorem states that for all
sufficiently large $j$, there exists a representation $\beta_j:
\Gamma\rightarrow {\rm PSL}_2({\bf C})$ so that $\beta_j(\Gamma)$ is a
minimally parabolic, geometrically finite uniformization of $M({\bf
p}_j,{\bf q}_j)$.  We refer to $\beta_j$ as the {\em hyperbolic Dehn
surgery representation of $\Gamma$ associated to $({\bf p}_j,{\bf
q}_j)$}. Moreover, the sequence $\{ \beta_j\}$ of representations of
$\Gamma$ converges to the identity representation of $\Gamma$. 

This form of the generalized hyperbolic Dehn surgery theorem is due to
Comar \cite{comar}.  There is also a version due to Bonahon and Otal
\cite{bonahon-otal}.  The original hyperbolic Dehn surgery theorem,
which applies in the case that the boundary of $M$ is equal to a union
of tori, is due to Thurston; the interested reader is directed to
Thurston's lecture notes \cite{thurston-notes} or to Benedetti and
Petronio \cite{benedetti-petronio}.  There is also the survey article
of Gromov \cite{gromov} discussing convergence and volumes of
hyperbolic $3$-manifolds. 

The main result of this Section is the observation that convex
co-compact subgroups of $\Gamma$ persist in the approximates
$\beta_j(\Gamma)$ given by the hyperbolic Dehn surgery theorem.

\begin{proposition} Let $M$ be a compact hyperbolizable $3$-manifold
and let $\{ T_1,\ldots, T_n\}$ be a non-empty collection of toroidal
components of $\partial M$.  Let $\Gamma$ be a minimally parabolic, 
geometrically finite uniformization of $M$. 

Let $(m_k,l_k)$ be a meridian-longitude system for $T_k$. Let $\{
({\bf p}_j, {\bf q}_j) = ((p^1_j,q^1_j), \ldots,(p^n_j,q^n_j)) \}$ be
a sequence of $n$-tuples of pairs of relatively prime integers such 
that, for each $k$, $\{ (p^k_j,q^k_j)\}$ converges to $\infty$ as
$j\rightarrow\infty$.   Let $\beta_j$ be the hyperbolic Dehn surgery
representation of $\Gamma$ associated to $({\bf p}_j,{\bf q}_j)$. 
 
Let $\Phi$ be a convex co-compact subgroup of $\Gamma$.  Then, for all
$j$ sufficiently large, $\beta_j(\Phi)$ is a convex co-compact
subgroup of $\beta_j(\Gamma)$.  Moreover, $\beta_j: \Phi\rightarrow
\beta_j(\Phi)$ is an isomorphism.
\label{persistence-prop}
\end{proposition}

\begin{proof} We begin the proof of Proposition \ref{persistence-prop}
with a definition.  Let $\varphi_1,\ldots, \varphi_j$ be a generating
set for $\Phi$.  An {\em $\varepsilon$-deformation} of $\Phi$ with
respect to this generating set is a representation $\psi$ of $\Phi$
into ${\rm PSL}_2({\bf C})$ such that $|\varphi -\psi(\varphi)|
<\varepsilon$ for $1\le k\le m$.   Here, for an element $m$ of ${\rm
PSL}_2({\bf C})$, we take $|m|$ to be the matrix norm of a matrix
representative of $m$. 

By the quasiconformal stability theorem of Marden \cite{marden}, there
exists a constant $\varepsilon_0 >0$ so that an
$\varepsilon$-deformation $\psi$ of $\Phi$ for $\varepsilon
<\varepsilon_0$ is induced by a quasiconformal homeomorphism of
$\overline{\bf C}$. In particular, $\psi$ is an isomorphism from 
$\Phi$ to $\psi(\Phi)$, and $\psi(\Phi)$ is convex co-compact.

Since $\{\beta_j\}$ converges to the identity representation of
$\Gamma$ into ${\rm PSL}_2({\bf C})$, we have that
$\{\beta_j(\varphi_k)\}$ converges to $\varphi_k$ for each $1\le k\le
m$.  In particular, for $\varepsilon =\frac{1}{2}\varepsilon_0$, there
exists $N >0$ so that the restriction $\beta_j |_\Phi$ of $\beta_j$ to
$\Phi$ is an $\varepsilon$-deformation of $\Phi$ for $j\ge N$, and we
are done. 
\end{proof}

We can recast this discussion slightly.  Suppose as above that $M$ is
a compact, hyperbolizable $3$-manifold, that $\{ T_1,\ldots, T_n\}$ is
a collection of toroidal components of $\partial M$, and that $\Gamma$
is a minimally parabolic, geometrically finite uniformization of $M$. 

Let $\cal P$ denote the set of all pairs $(p,q)$ of relatively prime
integers, together with $\infty$.  Identify the pair $(p,q)$ with
$p+q{\bf i}\in {\bf C}$, so that we can view $\cal P$ as a subset of
the Riemann sphere $\overline{\bf C}$ with its usual topology.  

For each $n$-tuple of pairs $({\bf p}, {\bf q}) =((p_1, q_1),
\ldots,(p_n, q_n))$ in ${\cal P}^n$, let $M({\bf p}, {\bf q})$ denote
the manifold  obtained by performing $(p_k, q_k)$ Dehn surgery along
$T_k$, relative to a chosen meridian-longitude system on $T_k$.  (If
the Dehn surgery coefficent in $\cal P$ for $T_k$ is $\infty$, then we
delete $T_k$ from $M$.)

Then, the generalized hyperbolic Dehn surgery theorem can be rephrased
as saying that there exists a neighborhood $U$ of $(\infty,\ldots,
\infty)$ in $\overline{\bf C}^n$ so that for every point $({\bf p},
{\bf q})$ in $U\cap {\cal P}^n$, there exists a representation
$\beta: \Gamma\rightarrow {\rm PSL}_2({\bf C})$ so that
$\beta(\Gamma)$ is a minimally parabolic, geometrically finite
uniformization of $M({\bf p},{\bf q})$.  Also, the Dehn surgery
representations converge to the identity representation of $\Gamma$
into ${\rm PSL}_2({\bf C})$ as $({\bf p}, {\bf q})$ converges to
$(\infty,\ldots, \infty)$. 

In the case that we are performing Dehn surgery along a single torus
component $T$ of $\partial M$, we can use this discussion to obtain a
variant of a Theorem of Wu \cite{wu}.  Building on work of Culler,
Gordon, Luecke, and Shalen \cite{cgls}, Wu proves the following.

\begin{theorem} (Wu \cite{wu}) Let $M$ be a compact orientable
$3$-manifold, let $T$ be a toroidal component of $\partial M$, and let
$S$ be a component of $\partial M -T$.  Suppose that $S$ is
incompressible in $M$ and that $S\cup T$ is an-annular.  Then, there
are at most three Dehn surgeries along $T$ so that $S$ is not
incompressible in the resulting $3$-manifold. 
\label{wu-theorem}
\end{theorem}

Using Proposition \ref{persistence-prop}, we can give a
non-quantitative proof  of Theorem \ref{wu-theorem} for
hyperbolizable $3$-manifolds.  As above, let $M$ be a compact
hyperbolizable $3$-manifold with a single toroidal boundary component
$T$, and let $\Gamma$ is a minimally parabolic, geometrically finite
uniformization of $M$.  Let $S$ be a incompressible component of
$\partial M -T$.  By the generalized hyperbolic Dehn surgery theorem,
all but finitely many of the $3$-manifolds obtained from $M$ by
performing Dehn surgery along $T$ are hyperbolizable.

Suppose that $S$ has genus $1$, so that $S$ is itself a toroidal
component of $\partial M$.   Then, $S$ is a component of the boundary
of the hyperbolizable $3$-manifold obtained by
performing hyperbolic Dehn surgery along $T$ for all but finitely many
Dehn surgery coefficients.  The observation that boundary tori of
hyperbolizable $3$-manifolds are necessarily incompressible completes 
this case.

Suppose that $S$ has genus at least two, and let $\Phi$ be a choice of
conjugacy class of $\pi_1(S)$ in $\Gamma =\pi_1(M)$.  The hypotheses
that $\Gamma$ is a minimally parabolic, geometrically finite
uniformization of $M$, that $S$ is incompressible, and that $S\cup T$
is an-annular together imply that $\Phi$ is a convex co-compact
quasifuchsian subgroup of $\Gamma$.  By Proposition
\ref{persistence-prop}, the Dehn surgery representation restricts to
an isomorphism of $\Phi$ with convex co-compact image for all but
finitely many Dehn surgeries along $T$.

Though we cannot obtain the quantitative information in this case that
Wu obtains in the general case, we are also not restricted to
subgroups of $\Gamma$ corresponding to incompressible surfaces in $M$,
but can instead show that every convex co-compact subgroup of $\Gamma$
persists in all but finitely many of the Dehn fillings along $T$.  In
particular, we have that a non-quantative version of Theorem
\ref{wu-theorem} for immersed surfaces holds.  

Also, we are not restricted to working in hyperbolizable $3$-manifolds
with a single cusp, but have a non-quantative version of Theorem
\ref{wu-theorem} for an embedded or immersed quasifuchsian surface in
an $n$-cusped finite volume hyperbolic $3$-manifold, and in fact in
any hyperbolizable $3$-manifold uniformized by a geometrically finite,
minimally parabolic Kleinian group, a class which includes many
infinite volume hyperbolic $3$-manifolds.

\section{Locally free groups}
\label{locally-free-groups}

A group $G$ is {\em locally free} if every finitely generated subgroup
of $G$ is free.  Every free group is locally free, as all subgroups of
free groups are free by Grushko's theorem.  Another standard example
of a locally free group is the additive group ${\bf Q}$ of rational
integers, in which every finitely generated subgroup is infinite
cyclic. 

There is an example of Maskit \cite{maskit-duke} of a Kleinian group
which is locally free but not free.  However, this example is
infinitely generated, and does not seem to reside inside a finitely
generated Kleinian group.

A fourth example, which we make heavy use of in this paper, is
contained in the group 
\[ G_0 =\langle a, b, t\: |\: t\cdot a\cdot t^{-1} =[b,a] \rangle, \]
where $[b,a] =b\cdot a\cdot b^{-1}\cdot a^{-1}$ is the commutator of
$b$ and $a$.  Namely, the subgroup
\[ H_0 =\langle t^m\cdot a\cdot t^{-m}, \: t^m\cdot b\cdot t^{-m}, \:
m\in {\bf Z} \rangle \]
of $G_0$ is locally free but not free.  We refer the Reader to Maskit
\cite{maskit-book}, Chapter VIII.E.9, and to Freedman and Freedman
\cite{freedman} for a more detailed discussion of this group, inluding
the proof that $H_0$ is locally free but not free.

We note here that Maskit realizes $G_0$ as a convex co-compact
Kleinian group $\Gamma_0$, and that the boundary of the associated
compact $3$-manifold $M_0 =({\bf H}^3\cup\Omega(\Gamma_0))/\Gamma_0$ 
is an incompressible surface of genus two.

We close this Section with the observation that the property of a
Kleinian group containing a subgroup which is locally free but not
free is preserved under commensurability.

\begin{proposition} Suppose that $\Phi_1$ and $\Phi_2$ are
commensurable Kleinian groups.  Then, $\Phi_1$ contains a subgroup
which is locally free but not free if and only if $\Phi_2$ contains a
subgroup which is locally free but not free.
\label{commensurable-wacky}
\end{proposition}

\begin{proof} Let $\Xi$ be a subgroup of $\Phi_1$ which is locally
free but not free, and let $\Theta =\Phi_1\cap\Phi_2$.  We first show
that $\Theta$ contains a subgroup which is locally free but not free,
namely the intersection $\Theta\cap\Xi$.

Every finitely generated subgroup of $\Theta\cap\Xi$ is a finitely
generated subgroup of $\Xi$, and so is free.  In particular, this
shows that $\Theta\cap\Xi$ is locally free.  

To see that $\Theta\cap\Xi$ is not free, we first note that since
$\Theta$ has finite index in $\Phi_1$, $\Theta\cap\Xi$ has finite
index in $\Xi$.  If $\Theta\cap\Xi$ is free, then $\Xi$ is a
torsion-free group containing a free subgroup of finite index, and so
by a Theorem of Stallings \cite{stallings} and Swan \cite{swan}, we
have that $\Xi$ is free, a contradiction.

Hence, $\Theta$ contains the subgroup $\Theta\cap\Xi$ which is locally
free but not free.  Since $\Theta$ is a subgroup of $\Phi_2$, we see
that $\Phi_2$ contains a subgroup which is locally free but not free. 
\end{proof} 

\section{A criterion for non-commensurability}
\label{non-commensurability-criterion}

The purpose of this Section is to discuss criteria for a collection
of hyperbolic $3$-manifolds to contain infinitely many commensurability
classes.  For a collection of finite volume hyperbolic $3$-manifolds,
one such criterion is given by an upper bound on the volume.

\begin{proposition} Let ${\cal M}$ be a collection of infinitely many
distinct finite volume hyperbolic $3$-manifolds.  Suppose there exists
some constant $K >0$ so that ${\rm vol}(N)\le K$ for all $N$ in ${\cal
M}$.  Then, no commensurability class in ${\cal M}$ contains
infinitely many manifolds.  In particular, ${\cal M}$ contains
infinitely many commensurability classes. 
\label{commensurability-prop}
\end{proposition}

\begin{proof} We begin the proof of Proposition
\ref{commensurability-prop} with a definition.  For a Kleinian group
$\Gamma$, define the {\em commensurability subgroup} ${\rm
comm}(\Gamma)$ of $\Gamma$ in ${\rm PSL}_2 ({\bf C})$ to be 
\[ {\rm comm}(\Gamma) =\{ g\in {\rm PSL}_2({\bf C})\: |\: \Gamma\mbox{
and }g\: \Gamma\: g^{-1}\mbox{ are commensurable.} \}. \]
Note that $\Gamma\subset {\rm comm}(\Gamma)$.  Moreover, if $\Gamma$
and $\Gamma'$ are commensurable, then ${\rm comm}(\Gamma) ={\rm
comm}(\Gamma')$.

Suppose there exist infinitely many manifolds $N_n ={\bf
H}^3/\Gamma_n$, $n\ge 1$, in $\cal M$ which are commensurable.  Since
$\Gamma_n$ and $\Gamma_m$ are commensurable, we have that ${\rm
comm}(\Gamma_n) ={\rm comm}(\Gamma_m)$ for all $n$, $m\ge 1$.  Set
$\Theta ={\rm comm}(\Gamma_n)$.  

Though we do not give a precise definition here, we note that there
exists a special class of Kleinian groups, the {\em arithmetic}
Kleinian groups, which roughly are Kleinian groups defined by number
theory.  For a detailed discussion of arithmeticity, we refer the
interested Reader to the forthcoming book of Maclachlan and Reid
\cite{maclachlan-reid}.  

For the purposes of this note, it suffices to make use of a major
result of Margulis \cite{margulis}, see also the discussion in Zimmer
\cite{zimmer}, which states that a finite co-volume Kleinian group
$\Gamma$ is arithmetic if and only if $\Gamma$ has infinite index in
${\rm comm}(\Gamma)$.  Note that arithmeticity of $\Gamma$ implies
that ${\rm comm}(\Gamma)$ is not discrete, and hence is dense in ${\rm
PSL}_2({\bf C})$.

We now apply a result of Borel \cite{borel}, which states that given a
constant $K >0$, there exist only finitely many arithmetic Kleinian
groups $\Phi$ with co-volume ${\rm vol}({\bf H}^3/\Phi)$ at most $K$.
In particular, only finitely many of the $\Gamma_n$ can be 
arithmetic, and so there must exist some $j\ge 1$ for which $\Gamma_j$
is not arithmetic.  In particular, $\Gamma_j$ is a subgroup of finite
index in $\Theta ={\rm comm}(\Gamma_j)$, and so $\Theta$ is discrete. 

We now make use of a fact about the set of volumes of hyperbolic
$3$-manifolds, namely that there is a minimum volume hyperbolic
$3$-manifold, and so there exists some $C>0$ so that ${\rm
vol}(N_n)\ge C$ for all $n$.  In particular, the index of $\Gamma_n$
in $\Theta$ is bounded above by 
\[ [\Theta: \Gamma_n] =\frac{{\rm vol}({\bf H}^3/\Theta)}{{\rm
vol}(N_n)}\le \frac{{\rm vol}({\bf H}^3/\Theta)}{C} \] 
for all $n\ge 1$.

In particular, infinitely many of the $\Gamma_n$ have the same index
in $\Theta$.  Since $\Theta$ is finitely generated, it has only
finitely many subgroups of a given finite index, and so infinitely
many of the $\Gamma_n$ are equal, a contradiction.  This contradiction
completes the proof.  
\end{proof}

At this point, there are several remarks to make.

\begin{remark} Note that Proposition \ref{commensurability-prop} also
holds, with the same hypotheses, for hyperbolic $2$-manifolds, as the
results of Margulis and Borel both apply to ${\rm PSL}_2({\bf R})$, as
does the lower bound of the volume (in this case, area) of the
quotient manifold (in this case, surface).
\label{2-dim-case} 
\end{remark}

\begin{remark} In the case that $N ={\bf H}^3/\Gamma$ is a finite
volume hyperbolic $3$-manifold with a single cusp, there is an
alternative proof of Proposition \ref{commensurability-prop} using
work of Long and Reid \cite{long-reid}, strengthening an earlier
observation of Hodgson. Namely, for each $D \in {\bf N}$, there are
(up to conjugacy) only finitely many discrete (not necessarily
faithful) representations $\rho$ of $\Gamma$ into ${\rm PSL}_2({\bf
C})$ for which the degree of the invariant trace field of
$\rho(\Gamma)$ has degree at most $D$.  Here, the {\em invariant trace
field} of $\rho(\Gamma)$ is the extension of ${\bf Q}$ generated by
the traces of the squares of the elements of $\rho(\Gamma)$. 

It is a result of Reid \cite{reid-trace}, that the invariant trace
field of a finite co-volume Kleinian group is an invariant of the
commensurability class.  Also, the Thurston hyperbolic Dehn surgery
theorem implies that the degree of the invariant trace field goes to
infinity under Dehn surgery, which then implies that Dehn surgery on a
singly cusped, finite volume hyperbolic $3$-manifold yields infinitely
many commensurability classes.

I would like to thank Joe Masters and Alan Reid for bringing this to
my attention. 
\label{second-remark}
\end{remark} 

\begin{remark} For hyperbolic $n$-manifolds with $n\ge 4$, Proposition
\ref{commensurability-prop} is a trivial consequence of a result of
Wang \cite{wang}, which states that for each $K >0$, there are only
finitely many isometry classes of hyperbolic $n$-manifolds of volume
at most $K$.
\end{remark}

One application of Proposition \ref{commensurability-prop} is to
the collections of hyperbolic $3$-manifolds obtained from a given finite
volume, cusped hyperbolic $3$-manifold $N$. 

\begin{corollary} Let $M$ be a compact hyperbolizable $3$-manifold and
suppose that $\partial M$ is a union of tori $\partial M
=T_1\cup\cdots\cup T_j$.  Let ${\cal M}(M)$ denote the set of all
finite volume hyperbolic $3$-manifolds obtained from $M$ by performing
Dehn surgery along some or all of the tori in $\partial M$.  Then, at
most finitely many of the manifolds in ${\cal M}(M)$ are
commensurable. 
\label{dehn-surgery-not-commensurable} 
\end{corollary}

\begin{proof} This Corollary follows immediately from Proposition
\ref{commensurability-prop} and the fact, due to Thurston, that the
volumes of the manifolds in ${\cal M}(M)$ are bounded by ${\rm
vol}(N)$; for a discussion of this fact, see for instance Gromov
\cite{gromov} or Benedetti and Petronio \cite{benedetti-petronio}. 
\end{proof}

Corollary \ref{dehn-surgery-not-commensurable} itself has a simple
Corollary, follows immediately from the observation that for each
$n\ge 0$, there exists a finite volume hyperbolic $3$-manifold with
$n+1$ cusps.   For $n =0$, Corollary
\ref{n-infinitely-many-classes} is due to Maclachlan and Reid
\cite{maclachlan-reid-comm} using arithmetic techniques. 

\begin{corollary} For $n\ge 0$, the collection ${\cal M}_n$ of all
finite volume hyperbolic $3$-manifolds with $n$ cusps contains
infinitely many commensurability classes.
\label{n-infinitely-many-classes}
\end{corollary}

There is also a variant of Proposition \ref{commensurability-prop} for
a collection of geometrically finite, infinite volume hyperbolic
$3$-manifolds.  First, recall that if $\Gamma$ is a geometrically
finite, infinite co-volume Kleinian group, then $\Omega(\Gamma)$ is
non-empty. 

\begin{proposition} Let ${\cal M}$ be a collection of infinitely many
distinct geometrically finite, infinite volume hyperbolic
$3$-manifolds.  Suppose there exists a constant $K >0$ so that ${\rm
area}(\Omega(\Gamma)/\Gamma)\le K$ for all $N ={\bf H}^3/\Gamma$ in
$\cal M$. Then, no commensurability class in ${\cal M}$ contains
infinitely many manifolds.  In particular, ${\cal M}$ contains
infinitely many commensurability classes. 
\label{commensurability-infinite-vol}
\end{proposition}

\begin{proof} Suppose there exist infinitely many elements $N_1$,
$N_2,\ldots, N_m,\ldots$ of $\cal M$ which are pairwise
commensurable.  Write $N_m ={\bf H}^3/\Gamma_m$, so that $\Gamma_m\cap
\Gamma_{m+1}$ has finite index in both $\Gamma_m$ and $\Gamma_{m+1}$.
Let ${\rm comm}(\Gamma_m)$ be the commensurability subgroup of
$\Gamma_m$.  Since $\Gamma_m$ and $\Gamma_{m+1}$ are commensurable for
each $m$, we have that ${\rm comm}(\Gamma_m) ={\rm comm}(\Gamma_p)$
for $m\ne p$.  Set $\Theta ={\rm comm}(\Gamma_m)$.

\medskip

We first show that $\Lambda(\Gamma_m) =\Lambda(\Theta)$.  Since
$\Gamma_m$ is a subgroup of $\Theta$, we immediately have that
$\Lambda(\Gamma_m) \subset\Lambda(\Theta)$.  

For the opposite inclusion, let $g\in \Theta ={\rm comm}(\Gamma_m)$ be
any element.  Since $\Gamma_m$ and $g\: \Gamma_m\: g^{-1}$ are
commensurable, we have that $\Lambda(\Gamma_m) = \Lambda(g\:
\Gamma_m\: g^{-1})$.  Since $\Lambda(g\: \Gamma_m\: g^{-1})
=g(\Lambda(\Gamma_m))$, we have that $g(\Lambda(\Gamma_m))
=\Lambda(\Gamma_m)$, and so $\Lambda(\Gamma_m)$ is a non-empty closed
subset of $\overline{\bf C}$ invariant under $\Theta$.  Since $\Theta$
contains $\Gamma_m$, no smaller non-empty  subset of $\overline{\bf
C}$ can be invariant under $\Theta$, and so $\Lambda(\Gamma_m)
=\Lambda(\Theta)$.

\medskip

There are now two cases.  In the case that $\Lambda(\Gamma_m)$ is not
a circle in $\overline{\bf C}$, it is a standard fact, see for
instance Sullivan \cite{sullivan}, that
\[ \Xi ={\rm stab}_{{\rm PSL}_2({\bf C})} (\Lambda(\Gamma_m)) \]
is a discrete subgroup of ${\rm PSL}_2({\bf C})$ with $\Lambda(\Xi)
=\Lambda(\Gamma_m)$.  

Since each $\Gamma_m$ is non-elementary, there is a canonical metric
of constant curvature $-1$ on $\Omega(\Gamma_m) =\Omega(\Theta)
=\Omega(\Xi)$ which is invariant under $\Xi$, and so descends to the
hyperbolic surfaces $\Omega(\Xi)/\Gamma_m$, $\Omega(\Xi)/\Theta$, and
$\Omega(\Xi)/\Xi$.  By the Ahlfors  finiteness theorem 
\cite{ahlfors}, $\Omega(\Xi)/\Gamma_m$ has finite area.  Since
$\Omega(\Xi)/\Gamma_m$ covers both $\Omega(\Xi)/\Theta$ and
$\Omega(\Xi)/\Xi$, both of these hyperbolic surfaces have finite area.
Since there is a lower bound on the area of a hyperbolic surface, this
implies that each $\Gamma_m$ has finite index in $\Xi$, and so $\Xi$
is finitely generated and discrete.  The index $[\Xi: \Gamma_m]$ is
bounded by the quotient of the areas ${\rm area}(\Omega(\Xi)/\Xi) /
{\rm area}(\Omega(\Gamma_m)/\Xi)$, and so as in the proof of
Proposition \ref{commensurability-prop} infinitely many of the
$\Gamma_m$ must be equal. 

In the case that $\Lambda(\Xi)$ is a circle in $\overline{\bf C}$, the
stabilizer of $\Lambda(\Xi)$ in ${\rm PSL}_2({\bf C})$ is conjugate to
${\rm PSL}_2({\bf R})$, and so we are done by Remark \ref{2-dim-case}.
\end{proof}

We close this Section by noting that the conclusion and proof of
Proposition \ref{commensurability-infinite-vol} also hold for
collections of finitely generated but geometrically infinite Kleinian
groups, as long as it is assumed that the domains of discontinuity of
the Kleinian groups in the collection are non-empty. 

\section{The construction}
\label{construction}

We are now ready to give the proof of Theorem \ref{infinitely-many}.

\begin{theorem} There exist infinitely many commensurability classes
of finite volume hyperbolic $3$-manifolds whose fundamental groups
contain a subgroup which is locally free but not free.
\label{infinitely-many}
\end{theorem}

\begin{proof} We begin by recalling the example of Maskit from
Section \ref{locally-free-groups} of a convex co-compact Kleinian group
$\Gamma_0$ which contains a subgroup which is locally free but not
free.  The boundary of the compact $3$-manifold $M_0 =({\bf
H}^3\cup\Omega(\Gamma_0))/\Gamma_0$ corresponding to $\Gamma_0$ is a
connected surface of genus two which is incompressible in $M_0$.

\medskip  

Let $M_1$ be a compact hyperbolizable acylindrical $3$-manifold whose
boundary $\partial M_1$ is the union of a torus $T$ and a surface $S$
of genus two.  Such manifolds are common.  To construct one, begin
with a compact hyperbolizable acylindrical $3$-manifold $M$ whose
boundary is a closed surface of genus two.  Such manifolds are known
to exist by work of Kojima and Miyamoto \cite{kojima-miyamoto}, see
also Paoluzzi and Zimmermann \cite{paoluzzi-zimmermann}.  Here, we are
using the fact that the condition that $M$ be acylindrical is
equivalent to the existence of a hyperbolic structure with totally
geodesic boundary on $M$.   

It is a result of Myers \cite{myers} that for each element $\gamma$ in
$\pi_1(M)$, there exists a simple closed curve $c$ in $M$ representing
$\gamma$ so that $M -N(c)$ is hyperbolizable, where $N(c)$ is a
regular neighborhood of $c$ in the interior ${\rm int}(M)$ of
$M$.  

We now consider how to choose $\gamma$ so that $M -N(c)$ is
acylindrical.  Let $(m,l)$ be a meridian-longitude system on the torus
$\partial N(c)$, as described in Section \ref{persistence}.  Choose
$\gamma$ so that $\gamma$ is not homotopic into $\partial M$.  Suppose
there exists an essential annulus $A$ in $M -N(c)$ so that one
component of $\partial A$ lies in $\partial M$ and the other component
of $\partial A$ lies in $T =\partial N(c)$. Write $\partial A =B\cup
C$, where $B$ lies in $\partial M$ and $C$ lies in $T$.

Express $C$ in the chosen meridian-longitude system, so that it is a $(p,q)$ curve
on $T$, where $p$, $q$ are relatively prime.  

If $C$ is a $(1,0)$ curve on $T =\partial N(c)$, then $C$ bounds a
disc $D$ in $N(c)$.  Then, the union $A\cup D$  is a properly embedded
disc in $M$ whose boundary $B$ is a homotopically non-trivial curve on
$\partial M$, which contradicts the incompressibility of $\partial M$.  

If $C$ is a $(p,q)$ curve on $T$ for $q\ne 0$, then $C$ is homotopic
to $c^q$ in $M$.  However, the acylindricality of $M$ implies that if
$c^q$ is homotopic into $\partial M$, then $c$ is homotopic into
$\partial M$, a contradition.  To see this, we use a result of Maskit
\cite{maskit-inter}, which states that if $\Gamma$ is a finitely
generated Kleinian group and if $\Delta$ and $\Delta'$ are two
components of $\Omega(\Gamma)$, then 
\[ \partial\Delta \cap\partial\Delta' =\Lambda({\rm
stab}_\Gamma(\Delta)\cap {\rm stab}_\Gamma(\Delta')). \]

So, if $c^q$ is homotopic into $\partial M$, then $\gamma^q$
stabilizes a component $\Delta$ of $\Omega(\Gamma)$, where $\Gamma$ is
a realization of $\pi_1(M)$ as a Kleinian group.  In particular, the
fixed points of $\gamma$ lie in $\partial \Delta$.  Since $M$ is
acylindrical, the components of $\Omega(\Delta)$ have disjoint
closures, see for instance Anderson \cite{anderson-ess}.  So, it
cannot be that $\gamma(\Delta)\ne \Delta$, as then the two components
$\Delta$ and $\gamma(\Delta)$ of $\Omega(\Gamma)$ would have disjoint
components with intersecting closures.  So, $\gamma(\Delta) =\Delta$,
and so $\gamma$ and hence $c$ are homotopic into $\partial M$.

The acylindricality of $M$ implies that there is no essential annulus
in $M -N(c)$ with both boundary components in $\partial M$, as then
would then give an essential annulus in $M$. Also, there can be no
essential annulus in $M -N(c)$ with both boundary components in
$\partial N(c)$.  Hence, there are no essential annuli in $M$.

\medskip

Let $h: \partial M_0\rightarrow S$ be an orientation-reversing
homeomorphism, and consider the $3$-manifold 
\[ M =M_0\cup_h M_1. \]
We now make use of Thurston's hperbolization theorem for
$3$-manifolds, see Morgan \cite{morgan}, which states that a compact,
orientable, irreducible, atoroidal $3$-manifold with non-empty boundary
is hyperbolizable.  As mentioned in the Introduction, this part of the
argument is standard and is included for the sake of completeness.  

Note that $M$ is compact, as both $M_0$ and $M_1$ are compact, and is
orientable, as $h$ is orientation-reversing.  Let $X$ be the image of
$\partial S =h(\partial M_0)$ in $M$. 

\medskip
\noindent
{\bf $M$ is irreducible: } Let $S$ be an embedded ${\bf S}^2$ in $M$.
First, isotope $S$ so that $S\cap X$ is the finite union of disjoint
simple closed curves.  Let $c$ be an innermost curve on $S$, meaning
that one of the components of $S -c$ contains no curve in $S\cap X$. 

Note that $c$ bounds a closed disc $D$ in $S$, namely the closure in
$S$ of the component of $S -c$ which contains no curve in $S\cap X$.
Since the interior of $D$ is disjoint from $X$, we have that $D$ is
contained in $M_k$, where either $k =0$ or $k = 1$.  Since $X$ is
incompressible in $M_k$ by assumption, $\partial D$ is a homotopically
trivial curve in $M_k$.  So, we can isotope $D$ into $S$ and thereby
get rid of $c$.  

Repeating this argument for each curve in $S\cap X$ in turn, working
from the innermost out, we can isotope $S$ into either $M_0$ or $M_1$.
Since both $M_k$ for $k=0$ and $k=1$ are irreducible, $S$ then bounds
a $3$-ball in $M_k$ and hence in $M$,  and so $M$ is irreducible as
well. 

\medskip
\noindent
{\bf $M$ is atoroidal: } This argument is very similar to the argument
that $M$ is irreducible.  Let $T$ be an incompressible torus in $M$,
and isotope $T$ so that $X\cap T$ is the finite union of disjoint
simple closed curves.  Again performing an innermost disc argument, we
can isotope away all the curves in $X\cap T$ which bound a disc in $T$.

Hence, we can assume that all the curves in $X\cap T$ are
homotopically non-trivial curves on the torus $T$.  Since both $X$ and
$T$ are embedded surfaces in $M$, the curves in $X\cap T$ are parallel
on $T$.  Moreover, since $X$ separates $M$, there must be an even
number of curves in $X\cap T$. If $X\cap T$ is empty, then $T$ is
contained in $M_k$ for either $k =0$ or $k =1$.  Since both $M_0$ and
$M_1$ are hyperbolizable, $T$ is then homotopic in $M_k$ into
$\partial M_k$. 

If $X\cap T$ is non-empty, then consider the closure $A$ of a
component of $T -(X\cap T)$ contained in $M_1$.  Since the curves in
$X\cap T$ are homotopically non-trivial on $T$, $A$ is an annulus.
Since $T$ is incompressible in $M$, the boundary curves are
homotopically non-trivial curves in $\partial M_1$, and so $A$ is an
incompressible annulus in $M_1$.  Since $M_1$ is acylindrical, we can
homotope $A$ into $\partial M_1$, and hence into $M_0$. 

Doing this for each component of $T -(X\cap T)$ which lies in $M_1$,
we can homotope all of $T$ into $M_0$.  Since $M_0$ is atoroidal, we
can homotope $T$ into a toroidal component of $\partial M_0$, which is
also a toroidal component of $\partial M$.  Hence, $M$ is atoroidal. 

\medskip

Hence, we have a hyperbolizable $3$-manifold, namely $M$, whose boundary
is a torus.  The Kleinian group $\Gamma$ uniformizing $M$ is
necessarily minimally parabolic and geometrically finite.  Also, since
$\Gamma$ contains a quasiconformal conjugate of $\Gamma_0$, we see
that $\Gamma$ contains a convex co-compact subgroup which in turn
contains a subgroup which is locally free but not free.  If we now
perform Dehn surgery along $\partial M$, we may combine Proposition
\ref{persistence-prop} and Proposition
\ref{commensurability-prop} to see that there exist
infinitely many commensurability classes of co-closed Kleinian groups
which contain a subgroup which is locally free but not free.
\end{proof}

Note that there is an extraordinary amount of flexibility in the
construction given in the proof of Theorem \ref{infinitely-many},
specifically in the choice of the $3$-manifold $M_1$ and in the choice
of the gluing map $h: \partial M_0\rightarrow S$.

\end{document}